\magnification=1200
\font\eightrm=cmr8
\font\eighti=cmti8
\font\eightrm=cmr8
\font\eighti=cmti8
\def\cqfd{\unskip\kern 6pt\penalty 500
\raise 0pt\hbox{\vrule\vbox to5pt{\hrule width 5pt
\vfill\hrule}\vrule}\par}
\overfullrule=0mm

{\parindent=0mm{\eightrm \footnote{}{2000 {\eighti Mathematics Subject Classification}.
Primary 37E30, Secondary 37E45. \obeylines
{\eighti Key words and phrases}. Brouwer Translation theorem, open annulus, fixed point, compactly connected component.}}} 
 \centerline{\bf ON THE STRUCTURE OF HOMEOMORPHISMS OF THE OPEN ANNULUS}

\bigskip
\centerline{LUCIEN GUILLOU }
\bigskip
\centerline{Dedicated to Jose Maria Montesinos on the occasion of his 65th birthday}
\bigskip
\bigskip

{\it Abstract}: Let $h$ be a without fixed point lift to the plane of a homeomorphism of the open annulus isotopic to the identity and without wandering point. We show that $h$ admits a $h$-invariant dense open set $O$ on which it is conjugate to a translation and we study the action of $h$ on the compactly connected components of the closed and without interior set ${\bf R}^2 \setminus O$.

\bigskip
\bigskip

\noindent{\bf 0. Introduction}.

\bigskip

0.1. In the paper [BCL] the authors consider homeomorphisms $H$ of the open annulus $S^1 \times {\bf R}$ isotopic to the identity and preserving the Lebesgue measure.  Given such a homeomorphism and a lift $h : {\bf R}^2 \rightarrow {\bf R}^2$ to the universal cover they show (in their proposition 3.1) that if the closure of the  rotation set of $h$ is contained in $]0, +\infty[$, then $h$ is conjugate to a translation. (Here the rotation set refers to a definition, adapted to this non compact situation, proposed by Le Calvez [LC] and using only recurrent points of $H$ in its construction).

They remark that this statement is sharp, and give an example of a measure preserving homeomorphism $H$ of $S^1 \times {\bf R}$ isotopic to the identity, such that, for some lift $h$ of $H$, the rotation set of $h$ is included in $]0, +\infty[$, but $h$ is not conjugate to a translation (see 0.2 below).

In the present note we wish to investigate the structure of such homeomorphims. More generally, we will consider a homeomorphism $H$ of the annulus $S^1 \times {\bf R}$ isotopic to the identity, without wandering point which admits a lift $h$ to ${\bf R}^2$ without fixed point. We will show that some of the features of example 0.2 are indeed preserved in that general situation.

\medskip

We will prove:

\medskip

\noindent{\it A) There exists an $h$-invariant dense open set homeomorphic to ${\bf R}^2$, $O \subset {\bf R}^2$, such that $h$ restricted to $O$ is conjugate to a translation.} (See paragraph 1).

\medskip

\noindent {\it B) Let $W = {\bf R}^2 \setminus O$ which is a closed subset with no interior in ${\bf R}^2$.}  We have:

\medskip

{\it - B1) No {\bf closed} compactly connected component (cf. 2.1 below) of $W$ is invariant under $h$.}  (Cf. Prop. 2.5).

\medskip

{\it - B2) We now assume that the compactly connected components of $W$ are {\bf closed} (or equivalently that the connected components of $W$ are compactly connected). Then for every such component $C$, if $X = \liminf \ h^n(C)$,

 either $X$ is empty (that is $h^n(C) \rightarrow \infty$, meaning that, for every compact $K \subset {\bf R}^2$ there exists an integer $n(K)$ such that $h^n (C) \bigcap K = \emptyset$ for $n \geq n(K)$) 

 or it is not empty and no point of $X$ is accessible from ${\bf R}^2 \setminus \overline { \bigcup _{n \in {\bf Z }} h^n(C)}$.}  (Cf. Prop. 2.12).

\bigskip

0.2. The Le Roux example [BCL, Appendix A]:

\medskip

We will describe the lift $h$ of this example to ${\bf R}^2$. Let $I_k$ be the vertical segment $\{ ({\displaystyle{1 \over 2k}}, y) \big \vert y \geq \vert k \vert \}$ and $A$ be  $\bigcup _{k \in {\bf Z} \setminus \{ 0 \}}I_k$ and let $W = \bigcup _{n \in {\bf Z}}T^n (A)$ where $T(x, y) = (x+1, y)$. Then ${\bf R}^2 \setminus W$ is homeomorphic to ${\bf R}^2$ and can be foliated by lines equivariantely with respect to $T$. The homeomorphism $h$ is choosen to act equivariantely, without fixed point,  preserving each line of the foliation and satisfying $h(I_k) = I_{k-1}$ for $k \not = 0, 1$. On each leaf of the foliation, $h$ is equivariantly conjugate to a translation hence $h$ preserves a measure without atoms and charging the open sets. On $S^1 \times {\bf R}$ seen as $S^2$ minus the two poles, $H$ preserves such a measure which is finite. That measure is nothing but the Lebesgue measure up to conjugation thanks to a classical result of Oxtoby and Ulam.

To see that h is not conjugate to a translation notice that the compact segment going from $x_0 = (-{1 \over 2}, 1)$ to its translate $T(x_0) = ({1 \over 2}, 1)$ 
 has to meet all its images by all iterates of $h$ since $W$ is $h$-invariant.

\medskip

We owe to P. Le Calvez the remark that this example can also be described without any reference to the Oxtoby-Ulam theorem. Consider the part of the phase space (which is homeomorphic to $S^1 \times {\bf R}$) of the free undamped pedulum above the upper separatrix: it is homeomorphic to $S^1 \times [0, +\infty [$. We now focus on the time 1 of the corresponding autonomous hamiltonian and on an orbit of this diffeomorphism on the separatrix. Folding each complementary interval of this orbit on the separatrix and identifying all points of the orbit and the equilibrium point of the separatrix to a single point, we get an example conjugate to the preceding one after deleting that single point.

\medskip

{\it Acknowledgement}. Some arguments of this paper can be traced back to an old article of T. Homma and H. Kinoshita [HK], which makes for a hard reading but a rewarding one. Many thanks to Patrice Le Calvez for a careful reading of a first version of this paper.

\bigskip

\bigskip

\noindent{\bf 1. Brouwer homeomorphisms}.

\bigskip

Homeomorphisms of the plane preserving orientation and without fixed point are called {\bf Brouwer homeomorphisms} (see [G1] for more on these).These homeomorphisms have only wandering points and more generally satisfy the following particular version  of Franks' lemma (in a reformulation due to Le Roux [LR1, Lemma 7]). Recall first that a subset $A$ of ${\bf R}^2$ is {\bf free} if $h(A) \bigcap A = \emptyset$.

\proclaim 1.1 Lemma. Let $U$ and $V$ be two free connected open sets.Then the subset of integers such that $h^n(U) \bigcap V \not = \emptyset$ is an interval of ${\bf Z}$.

{\it Proof}: The usual formulation of this lemma concerns the case where $U$ and $V$ are open discs. To prove the present lemma from this case, suppose there exists $k<n<m$ such that $h^k(U) \bigcap V \not = \emptyset, h^n(U) \bigcap V = \emptyset$ and $h^m(U) \bigcap V \not = \emptyset$. Let $u_1 \in U$ such that $v_1 = h^k(u_1) \in V$ and $u_2 \in U$ such that $v_2 = h^m(u_2) \in V$ and let $D$ and $D'$ be discs in $U$ and $V$ respectively such that $u_1, u_2 \in D$ and $v_1, v_2 \in D'$. Then $h^k(D) \bigcap D' \not = \emptyset, h^n(D) \bigcap D' = \emptyset$ and $h^m(D) \bigcap D' \not = \emptyset$ in contradiction to Franks' lemma.

\medskip

A {\bf Brouwer line} for a Brouwer homeomorphism $h$ is a properly embedded free line $l$ such that $l$ separates $h^{-1}(l)$ and $h(l)$.We will start with  the following result from [G2] .

\proclaim 1.2 Theorem. Let $H : S^1 \times {\bf R} \rightarrow S^1 \times {\bf R}$ be a homeomorphism isotopic to the identity such that : 

{\it - $H$ admits a fixed point free lift $h : {\bf R}^2 \rightarrow {\bf R}^2$.

- $H$ does not have any wandering point.

Then there exists a properly embedded line in $S^1 \times {\bf R}$ joining one end of the annulus to the other which lifts in ${\bf R}^2$ to a Brouwer line.}

\medskip

Notice that such a Brouwer line projects properly and onto on $\{ 0 \} \times {\bf R}$ (and also, a properly embedded line in ${\bf R}^2$ which projects properly and onto on $\{ 0 \} \times {\bf R}$ is a Brouwer line if it is free, that is, the requirement that $l$ separates $h^{-1}(l)$ and $h(l)$ is automatically satisfied).

\bigskip

Given any Brouwer line $l$, if we let $U$ be the open region between $l$ and $h(l)$, then the set $O = \bigcup _{n \in {\bf Z}}h^n (\hbox{Cl} U)$ is homeomorphic to ${\bf R}^2$ and the restriction of $h$ to $O$ is conjugate to a translation.Therefore to prove statement A of the introduction, it is enough to prove that if the Brouwer homeomorphism $h$ is a lift of a homeomorphism $H$ of the open annulus without wandering point, then ${\bf R}^2 \setminus O$ has no interior for a convenient choice of Brouwer line $l$. To this end, we choose a Brouwer line $l$ as given by Theorem 1.2 that we orient so that $l$ induces by projection the usual orientation on $\{ 0 \} \times {\bf R}$. The following Lemma is then enough to conclude the proof of statement A (this lemma is an extension of the  lemma in Winkelnkemper [W]).

\bigskip

\proclaim 1.3 Lemma. Let $B_n$ (resp. $B'_n$) be the component of ${\bf R}^2 \setminus
h^n (l)$ to the right (resp. to the left) of $h^n (l)$. Then the closed
$h$-invariant set $ W=\bigcap ^{+\infty}_{n=-\infty}B_n $ (resp. $ W'=\bigcap
^{+\infty}_{n=-\infty}B'_n $) has no interior. 

{\it Proof}:  
Exchanging $h$ and $h^{-1}$ if necessary, we can suppose $h(l)$ on the right of $l$. Suppose $U\subset W$ is an open subset which we
can choose small enough to be free and projecting homeomorphically on $S^1\times {\bf R}$; since $U \subset W$, $h^{-n}(U)$ lies on the right of $l$ for all $n \geq 0$. Given
the properties of $l$, there is a $m>0$ such that $U$ lies on the left of $T^m(l)$, then $h^{-n}(U)$
lies on the left of $h^{-n}(T^m(l))=T^m(h^{-n}(l))$ which is on the left of $T^m(l)$ for $n>0$. So that
all $h^{-n}(U)$, $n \geq 0$, lie on the left of $T^m (l)$ and on the right of $l$. There are only a finite
number of translates of $ U$ between $l$ and $T^m (l)$, say $U=U_1 , U_2 ,\ldots , U_k$ and each
one is wandering. Since $H$ has no wandering point on $S^1\times {\bf R}$, there exists $n_1>0$ such
that $h^{-n_1}(U_1)$ meets some $U_i$, say $U_{j(1)}$. Let $V_1= h^{-n_1}(U_1) \bigcap U_{j(1)}$.
There exists also $n_2>0$ such that $h^{-n_2}(V_1)$  meets one of
its translates $V_{j(2)} \subset U_{i_2}$. Let $V_2=h^{-n_2}(V_1) \bigcap V_{j(2)}$. 
Continuing in that way we find a sequence $V_1, V_2, \ldots$ of non empty sets each $V_i$ being
contained in some $U_{j(i)}$, $1 \le j(i) \le k$. We must have $j(i) = j(i')$ for some $i$ and $i'$,
$i<i'$. Then, since $V_{i'} \subset h^{-p}(V_i)$ for $p =  n_{i+1} + \ldots + n_{i'}$, we have
$U_{j(i')} \bigcap h^{-p}(U_{j(i)}) \not = \emptyset$ contradicting the freeness of $U_{j(i)}$.

\bigskip

\noindent{\bf 2. Compactly connected components}

\bigskip

In this paragraph we consider any Brouwer homeomorphism $h$ and an associated oriented Brouwer line
$l$ such that $ W=\bigcap ^{+\infty}_{n=-\infty}B_n $ and $ W'=\bigcap
^{+\infty}_{n=-\infty}B'_n $ have no interior (where as above $B_n$ (resp. $B'_n$) is the
component of ${\bf R}^2 \setminus h^n (l)$ to the right (resp. to the left) of $h^n (l))$.

  Notice that the sets $W$ and $W'$ are disjoints, that  the invariant set $O = {\bf R}^2
\setminus (W�\bigcup W')$ is homeomorphic to ${\bf R}^2$ and that on this set $h$ is conjugate to a
translation. Similar considerations can be applied to each one of $W$ and $W'$ and we will only
describe those pertaining to $W$.

The set $W$ is generally not connected. It is also non-compact (since  it is invariant and points are wandering under $h$) and we will have to consider its compactly connected components. Let us recall (see [Moore, page 76] and also [LR2, D\'efinition 9.1])

\proclaim 2.1 Definition. A space $Z$ is compactly connected if any two points in $Z$ are
contained in a subcontinuum of $Z$. Distinct maximal compactly connected subsets of a space $X$ are disjoint
and are called the compactly connected components of $X$; these
components fill in $X$. Notice that these compactly connected components can be non closed.

\proclaim 2.2 Lemma. The compactly connected components of $W$ are unbounded.

{\it Proof}: We work in the Alexandroff compactification of ${\bf R}^2$, that is ${\bf R}^2 \bigcup \{ \infty \} \cong S^2$. First, $W \bigcup \{ \infty \}$ is compact and connected as the decreasing intersection of the compact connected $B_n \bigcup \{ \infty \}$. Suppose now that $W$ admits a compactly connected component $C$ contained in some open ball $B(O, R)$. Then $C$ is connected and compact so is a connected compact component of $W$. As such, it is the intersection of the open and closed subsets of $W$ which contains $C$ [B, II \S 4.4], and there exists an open and closed neighborhood of $C$ inside $W \bigcap B(O, R)$. But this contradicts the connectivity of $W \bigcup \{ \infty \}$.

\medskip

Let us call $C$ a 
{\bf closed} compactly connected component of $W$ and $p$ an accessible point
of $C$ from ${\bf R}^2 \setminus C$ : $p$ is the extremity of an arc $\gamma $ such that $\gamma \setminus \{ p \} \subset {\bf R}^2 \setminus C$. We can suppose that 
$\gamma $ is
a free simple arc. Each $h^n (l)$ has to meet $\gamma $ and  $h(\gamma )$ for $n$ larger than some
$n_0$ which we can suppose to be $-1$, replacing $l$ by $h^{n_0 +1} (l)$ if necessary. Let $p_n$
denote the last point of $h^n (l)$ on $\gamma $ as we move towards $p$. Then the arc $\gamma
_n = p_n p$ on $\gamma $ is disjoint from all $h^i (l), i\leq n$ except for $p_n \in h^n (l)$.

Let $q_0 = h(p_{-1})$ and $\alpha_0$ be the subarc $p_0q_0$ of $l$. Since
 ${\bf R}^2 \setminus (W�\bigcup W')$ is simply connected (even homeomorphic to ${\bf R}^2$), it is divided by the arc \hfill\break
$\gamma _0 \bigcup \alpha_0 \bigcup h(\gamma _{-1})$ into two domains and we call
$\Omega $ the one which does not contain $h^{-1} (l)$.

\proclaim 2.3 Proposition. The domain $\Omega $ is free.

{\it Proof}: Suppose there exist $x \in \Omega \bigcap h(\Omega )$ and let $\beta$ be an arc from $a$ to $h^{-1}(x)$ with $a \in \hbox {int} p_0p$ and $\beta \setminus \{ a \} \subset \Omega$. Since $h$ preserves orientation, $h(y) \notin \Omega $ for $y$ close to $a$ on $\beta$. As $h(\beta) \bigcap h(p_0p) = h(a)$ and $h(\beta) \bigcap \alpha _0 = \emptyset$ (since $h(\beta)$ is on the right of of $h(l)$ and so, on the right of $l$ which contains $\alpha _0$), there exist some $b \in \beta$  such that the subarc $h(ab)$ of $h(\beta)$ joins $h(p_0p)$ to $p_0p$ inside ${\bf R}^2 \setminus (W \bigcup W' \bigcup \Omega )$ and the Jordan curve $\alpha _0 \bigcup q_0h(a) \bigcup h(ab) \bigcup h(b)p_0$ contains the whole Brouwer line $l$ or $C$ (according to $p_0$ or $p$ is contained inside that Jordan curve) which is absurd since these sets are unbounded.

\proclaim 2.4 Proposition. The closed compactly connected component $C$ cannot be $h$-invariant.

{\it Proof}: Assume by contradiction that $h(C) = C$ and let then $\widetilde K \subset C$ be a continuum containing $p$ and $h(p)$. Then $\Omega$ is bounded and being simply connected has a boundary Fr$\Omega$ which is connected and separating the plane. We first show
 $K = \overline { \Omega } \bigcap \widetilde K$ is compact and connected. It is enough to show that Fr$\Omega \bigcap C$ is connected for then, if $K = (\hbox{Fr}\Omega \bigcap C) \bigcap \widetilde K$ is not connected then (Fr$\Omega \bigcap C) \bigcup \widetilde K \subset C$ separates the plane which contradicts the fact that $C$ has no interior and does not separate. Let us note $\delta = \gamma _0 \bigcup \alpha _0 \bigcup h(\gamma _{-1})$ so that Fr$\Omega \bigcap C = \hbox{Fr}\Omega \setminus (\delta \setminus \{p, h(p) \}$. If this last set is not connected, it has either three components or more, and then Fr$\Omega $ is not connected or two components, containing $p$ and $h(p)$ respectively, which do not disconnect the plane and then Fr$\Omega$ does not disconnect.

\medskip

Therefore $\Sigma = \overline {\bigcup _{n \in {\bf Z}}h^n(K )} \subset W$ is a closed connected set which is invariant under $h$ and therefore non compact. As $W$ does not separate ${\bf R}^2$ and has no interior, the same is true of $\Sigma$ and ${\bf R}^2 \setminus \Sigma$ is homeomorphic to ${\bf R}^2$. The proper line $l$ separates  ${\bf R}^2 \setminus \Sigma$ into two regions homeomorphic to ${\bf R}^2$ and  we name $R$ the one between $l$ and $\Sigma$. The region $R$ itself is cut by the arc $p_0p$ 
 into two regions $A$ and $B$  where we call $A$ the one containing $\Omega$ and $B$ the one containing $h^{-1}(\Omega) \bigcap R$. By definition $p_0 p$ is on the frontier of $A$ and $B$. Notice that $A$ (and $B$) are non compact since we can follow $l$ to infinity in one direction or the other staying in $A$ (or $B$). Note that $A$ contains $h^k(\Omega), k \geq 0$ and $B$ contains $h^{-k}(\Omega) \bigcap R, k \geq 1$.

\medskip

\proclaim 2.5 Lemma. Fr$A$$\bigcap$Fr$B$$\bigcap \Sigma$ is non compact.

{\it Proof}: Let $\Sigma _A$ (resp. $\Sigma _B$) be the set of points of $\Sigma$ which admit a neighborhood contained in $A \bigcup \Sigma$ (resp. $B \bigcup \Sigma$). The sets $A \bigcup \Sigma _A$ and $B \bigcup \Sigma _B$ are disjoint and open, therefore their complement in $R \bigcup \Sigma \bigcup \setminus (p_0p \setminus \{ p \})$ (which is the set of points of $\Sigma $ for which every neighborhood meets $A$ and $B$, that is  Fr$A$$\bigcap$Fr$B$$\bigcap \Sigma$) separates $R \bigcup \Sigma \setminus (p_0p \setminus \{ p \})$ and $R \bigcup \Sigma \setminus (p_0p \setminus \{ p \})$ can be written as the disjoint union $(A \bigcup \Sigma _A) \coprod (B \bigcup \Sigma _B) \coprod ( \hbox{Fr}A \bigcap \hbox{Fr}B \bigcap \Sigma)$.

On the other hand, if Fr$A$$\bigcap$Fr$B$$\bigcap \Sigma$ was compact in ${\bf R}^2$ or equivalently in $R \bigcup \Sigma$ (which is homeomorphic to ${\bf R}^2$), thinking of $l$ as a straight line and of $p_0p$ as a segment orthogonal to $l$ (as it is legitimate by Schoenflies theorem), one can find a large rectangle in $R \bigcup \Sigma$ with a side parallel to $l$ , containing Fr$A$$\bigcap$Fr$B$$\bigcap \Sigma$  and whose boundary cuts $p_0p$ transversaly in a single point. The boundary of this rectangle joins a point of $A$ near $p_0p$ to a point of $B$ near $p_0p$ in contradiction to the above decomposition of  $R \bigcup \Sigma \setminus (p_0p \setminus \{ p \})$.

\medskip

Given Lemma 2.5, let us choose some point $x$ in  Fr$A$$\bigcap$Fr$B$$\bigcap \Sigma$ and outside $K$. Then $x \notin \overline {\Omega}$ and we choose an open euclidean ball $2U \subset {R \bigcup \Sigma}$ centered at $x$ free and disjoint of $\overline {\Omega}$. ($U$ will denote the ball of radius one half the one of $2U$). As $x$ belongs to $\Sigma$, $U$ meets some $h^m(K)$ and so some $h^m(\Omega)$ and (exchanging $h$ and $h^{-1}$ if necessary) we can suppose $m>0$ and therefore that $h^m (\Omega) \subset A$. Since $U$ meets $B$, we want to show that $2U$ meets some $h^{-n}(\Omega)$, for some $n>0$, for then $2U$ and $\Omega$ will give a contradiction to Lemma 1.1.

To that end, let us choose on Fr$U$ two arcs, one on Fr$U \bigcap A$ and the other on Fr$U \bigcap B$ (these exist since $U$ meets $A$ and $B$ which are connected non compact) and choose an arc $\alpha _0$ inside $R \setminus (\Sigma \bigcup U)$ joining these two arcs and meeting transversally  $p_0p$ into a single point. Complete $\alpha _0$ by a sub-arc $\alpha _1$ of Fr$U$. This gives a Jordan curve $\alpha$ inside $R \bigcup \Sigma$ which contains $p$ in its interior. Since points are wandering there exists $N>0$ such that $h^{-N}(p) \in \Sigma$ belongs to the exterior of $\alpha$.

Now, if $U$ does not meet any $h^{-k}(\Omega), k>0$, the connected set $\hat K = \bigcup_{i=1}^Nh^{-i}(K)$ either joins $p$ inside $\alpha$ to $h^{-N}(p)$ outside $\alpha$ without meeting $\alpha$ (in contradiction to the Jordan curve theorem, or it meets $\alpha _1$ ($\hat K$, contained in $\Sigma$, does not meet $\alpha _0$) and then $\overline {U}$ meets some $h^{-i}(K) \subset \hat K$ and so $2U$ meets some $h^{-i}(\Omega))$ and we are done. This concludes the proof of Proposition 2.4.

\bigskip

\proclaim 2.6 Corollary. $h^n (C) \bigcap C = \emptyset$ for all $n \in {\bf Z} \setminus \{ 0 \}$.

{\it Proof}:  If $h^n(C) \bigcap C \not = \emptyset$ then $h^n(C) = C$ in contradiction to 2.4 applied to $h^n$ which has the same $W$ as $h$.

\medskip

Recall that given given a sequence $\{X_n\}$ of subspaces of a topological space $Z$, a point $x \in Z$ belongs to $\liminf X_n$ if every neighborhood of $x$ meets $X_n$ for an infinite number of $n$ and to $\limsup X_n$ if every neighborhood of $x$ meets $X_n$ for all but a finite number of $n$.

\medskip

 We will now suppose that $X = \liminf h^n(C)$ is not empty.  It is then a closed and non compact subset of $W$  (since it is $h$-invariant). We aim to Proposition 2.12 below. Our first step is :

\medskip

\proclaim 2.8. Proposition. The set $X$ is also $\limsup h^n(C)$. That is, every open set $U$ which meets an infinite number of $h^n(C)$, meets $h^n(C)$ for all $n$ greater than some $n_0 = n_0 (U)$.

\medskip

\noindent {\it Remark}: This Proposition answers a question of F. Le Roux [LR2, footnote 7]

\medskip

{\it Proof}: We will use repeatedly the following immediate consequence of a result of Le Roux [LR2, Lemme 9.3], we repeat the proof here for completeness.

\medskip

\proclaim 2.9. Proposition. $X \bigcap h^n(C) = \emptyset$ for all $n \in {\bf Z}$.

{\it Proof}: Since $X$ is $h$-invariant, it is enough to show that $X \bigcap C = \emptyset$. Let us suppose $X \bigcap C \not = \emptyset$, and let $U$ be a free neighborhood of $x \in X \bigcap C$ such that $ {\overline U} \bigcap h(C) = \emptyset$. As $x \in X$, there exists $n>1$ so that $U \bigcap h^n(C) \not = \emptyset$. Let $y \in C$ such that $h^n(y) \in U$. There exists a continuum $K \subset C$ which contains $x$ and $y$. Since $h(C)$ (as $C$) is free, we can find a free connected neighborhood $V$ of $h(K) \subset h(C)$ such that $U \bigcap V = \emptyset$. But $x \in U \bigcap h^{-1}(V)$ and $h^n(y) \in U \bigcap h^{n-1}(V)$ so that $U$ and $V$ contradict Lemma 1.1.

\medskip

Let $V$ be a free open disc and $D$ a component of $V \setminus \overline{\bigcup _{n \in {\bf Z}} h^n(C)}$.

\medskip

\proclaim 2.10. Lemma. If Fr$D$ meets $h^n(C)$ and $h^m(C)$, then $\vert n-m \vert \leq 1$ and Fr$D$ cannot meet $X$ if it meets some $h^n(C)$.

{\it Proof}: To prove the first assertion, note that since  $X \bigcap h^n(C) = \emptyset$ for all $n$, given $x \in h^n(C)$ there exists a disc neighborhood $U$ of $x$ which does not meet any other $h^p(C)$ and a ray from $x$ to some point in $D \bigcap U$ leads to an accessible point of $h^n(C)$ from $D$. So let us suppose $\vert n-m \vert >1$ and let $\alpha$ be an arc from $a \in h^n(C)$ to $b \in h^m(C)$ such that $\alpha \setminus \{a, b \} \subset D$ and let $K$ be a continuum in $h^n(C)$ containing $a$ and $h^{n-m}(b)$. We assert that $K \bigcup \alpha$ is free. Indeed, $K$ is free as a subset of $h^n(C)$, $\alpha$ is free as $V$ is free and $h(K) \bigcap \alpha = \emptyset = h^{-1}(K) \bigcap \alpha$ since $n \pm 1 \not = m$. But $b \in h^{m-n}(K \bigcup \alpha) \bigcup (K \bigcup \alpha)$, and a small enough neighborhood of $K \bigcup \alpha$ will contradict Lemma 1.1 if $\vert n-m \vert >1$.

Let us suppose now that $X$  meets Fr$D$ and some $h^n(C)$ and let again $U$ be a disc neighborhood of some point $x \in \hbox{Fr}D \bigcap X$ small enough so that $U \bigcap h^k(C) = \emptyset$ if $\vert k \vert \leq \vert n \vert +1$. A ray issued from $x$ will either give an accessible point of some $h^m(C), \vert m \vert > \vert n \vert +1$ from $D$, but this is impossible according to the first part of the proof, or an accessible point of $X$ from $D$. In that case, let $\alpha$ be an arc from some point $a \in h^n (C)$ to $b \in X$ with $\alpha \setminus \{ a, b \} \subset D$ and let $U'$ be a free neighborhood of $b$ such that $U' \bigcap h^k (C) = \emptyset$, for $\vert k \vert \leq \vert n \vert + 1$ and such that $U' \bigcap h^{\pm 1}(\alpha) = \emptyset$. The arc $\alpha$ can be extended to an arc $\tilde {\alpha} \subset \alpha \bigcup U$ which joins $a \in h^n(C)$ to some $\tilde {b} \in h^m(C)$, $\vert m \vert > \vert n \vert +1$. If $K \subset h^n(C)$ is a continuum containing $a$ and $h^{n-m}(\tilde b)$, then $K \bigcup \tilde {\alpha}$ is a free continuum such that $\tilde b \in h^{m-n}(K \bigcup \tilde {\alpha}) \bigcap K \bigcup \tilde {\alpha}$ and a free neigborhood of this continuum gives a contradiction to Lemma 1.1.

\medskip

We now return to the proof of Proposition 2.8. Let $V$ be a free neighborhood of $x \in \hbox {liminf}\  h^n(C)$. There exist $m$ and $n>m+1$ such that $V$ meets $h^n(C)$ and $h^m(C)$. Let $\alpha$ be an arc in $V$ going from $a_m \in h^m(C)$ to $a_n \in h^n(C)$ disjoint from $h^m(C)$ and $h^m(C)$ except for its extremities. Let $D$ be the the component of $V \setminus \overline {\bigcup _{n \in {\bf Z}}h^n(C)}$ which meets $\alpha$ and has $a_n$ on its frontier. By 2.11, Fr$D$ meets $h^{n+1}(C)$ or $h^{n-1}(C)$. In the first case, let $a_{n+1}$ be the last point of $h^{n+1}(C)$ seen on $\alpha$ when going from $a_n$ to $a_m$. If $D'$ is the component of  $V \setminus \overline {\bigcup _{n \in {\bf Z}}h^n(C)}$  which meets the subarc $a_m a_{n+1}$ of $\alpha$ and has $a_{n+1}$ on its frontier, then Fr$D'$ does not meet $h^n(C)$ by construction of $\alpha$ and therefore, according to 2.10, meets $h^{n+2}(C)$. Iterating this process we see that $\alpha$ meets all $h^k(C), k \geq n$. In the other case, the same reasonning shows that $\alpha$ meets all the $h^k(C)$ for $m \leq k \leq n$. As $\alpha$ meets an infinite number of $h^k(C)$ we conclude in either case that $V$ meets all $h^k(C)$ for $k$ large enough and therefore $x \in  \limsup  h^n(C)$.

\medskip

\noindent {\bf  2.11. Assumption: We assume for the rest of this paper that the compactly connected components of $W$ (in fact, we will only consider those of $X$) are closed.}

\medskip

\proclaim 2.12. Proposition. No point of $X$ is accessible from  ${\bf R}^2 \setminus \overline {\bigcup _{n \in {\bf Z}}h^n(C)}$.

{\it Proof}: We begin with a lemma :

\medskip

\proclaim 2.13. Lemma. There is no free arc $\alpha$ joining $C$ to $X$ contained in ${\bf R}^2 \setminus \overline {\bigcup _{n \in {\bf Z}}h^n(C)}$ except for its extremities.

{\it Proof of 2.13}: Let $\alpha$ join $p \in C$ to $q \in X$ and consider a neighborhood $D$ of $q$ such that $\alpha \bigcup D$ is still free and $D \bigcap h(C) = \emptyset = D \bigcap h^{-1}(C)$ (recall that $X$ is disjoint from $h(C)$ and $h^{-1}(C)$ by proposition 2.9). Then $\alpha \bigcup D$ contains a point $h^n(p')$ for some $n>1$ and some $p' \in C$. Let $K \subset C$ be a continuum containing $p$ and $p'$ and consider the continuum $L = K \bigcup \alpha \bigcup D$. It is free but $h^n(p') \in h^n(L) \bigcap L$ and a small enough neighborhood of $L$ gives a contradiction to Lemma 1.1.

\medskip

 At this point we will finish the proof of 2.12 following the lines of the proof of a similar result (with $C$ replaced by a disc) in [LR2, Proposition 5.5].

\medskip

Let us suppose there exist a point $q$ of $X$ accessible from ${\bf R}^2 \setminus \overline {\bigcup _{n \in {\bf Z}}h^n(C)}$ by some arc $\alpha$ and let $Z$ be the connected component of $X$ which contains $q$. A point $x$ of ${\bf R}^2 \setminus \overline {\bigcup _{n \in {\bf Z}}h^n(C)}$ will be called a {\it neighborhood point} of $Z$ if there exists a free closed euclidean disc $D$ with center $x$ such that int$D \bigcap Z \not = \emptyset$. The set of all such points is an open set $V$

A point of $x \in V$ will be said of type $C$ if there is some euclidean disc $D$ of center $x$ as in the previous definition and an arc in $D$ from $x$ to $Z$ which meets some $h^n(C)$ and of type $Z$ if there exists such a disc $D$ and an arc in $D$ from $x$ to $Z$  which does not meet any $h^n(C)$. It follows from Lemma 2.13 that this type is well defined.

We show that all points of $V$ are of type $C$. Indeed, it is easily verified that the type is locally constant on $V$ and so is constant on every connected component of $V$. But $V \bigcup Z$  and ${\bf R}^2 \setminus Z$ are connected and therefore their intersection $V$ also as follows from the Mayer-Vietoris sequence of the pair $({\bf R}^2 \setminus Z, V \bigcup Z)$. Furthermore, since $Z \subset X$, certainly $V$ meets some $h^n(C)$ and all points of $V$ are of type $C$.

Now, if the point $x$ on the arc $\alpha$ is close enough to $q$, the subarc $xq$ of $\alpha$ is contained in a free euclidean disc which meets $Z$, and, $x$ being of type $C$, meets some $h^n(C)$. Contradiction.

 \bigskip \bigskip
{\eightrm
\centerline {REFERENCES}\bigskip
{ \parindent=5mm  
\item{[{ Br}]}   L.E.J. BROUWER, {\eighti Beweis des ebenen Translationssatzes}, Math.
Ann., 72 (1912), 37-54.

\item{[{ BCL}]} F. B\'EGUIN, S. CROVISIER and F. LEROUX, {\eighti Pseudo-rotations of the open annulus}, 
Bull. Braz. Math. Soc. 37 (2006), 275-306.

\item{[{ B}]} N. BOURBAKI,  Topologie g\'en\'erale, Hermann, Paris (1971)

\item{[{ G1}]}  L. GUILLOU, {\eighti Th\'eor\`eme de translation plane de Brouwer et
g\'en\'eralisations du th\'eor\`eme de Poincar\'e-Birkhoff}, Topology, 33 (1994), 331-351. 

\item{[{ G2}]} L. GUILLOU, {\eighti Free lines for homeomorphisms of the open annulus}, Trans. AMS, 360 (2008), 2191-2204.

\item{[{ HK}]} T. HOMMA, H. TERASAKA, {\eighti On the structure of the plane translation of Brouwer},
Osaka Math. J. 5 (1953) 233-266.

\item{[{ LC}]} P. LE CALVEZ, {\eighti Rotation numbers in the infinite annulus}, Proc. Amer. Math. Soc. 129 (2001), 3221-3230.

\item{[{ LR1}]} F. LEROUX, {\eighti Bounded recurrent sets for planar homeomorphisms}, Ergodic theory Dynam. Systems 19 (1999) 1085-1091.

\item{[{ LR2}]} F. LEROUX, {\eighti Structure des hom\'eomorphismes de Brouwer}, Geometry and Topology 9 (2005) 1689-1774.


\item{[{ M}]}  R.L. MOORE,   Foundations of Point Set Theory, AMS Colloquium Publications volume XIII, 1962

\item{[{ W}]} H. WINKELNKEMPER, {\eighti A generalisation of the Poincar\'e-Birkhoff
theorem}, Proc. AMS, 102 (1988), 1028-1030.   

\par}

\bigskip

Lucien Guillou 
 
Universit\'e Grenoble 1, Institut Fourier B.P. 74, Saint-Martin-d'H\`eres
38402  France 

E-mail: lguillou@ujf-grenoble.fr

\end